\documentclass[a4paper,11pt,reqno]{amsart}

\usepackage{amssymb,amsthm,amsmath}
\newtheorem{secc}{}
                         \newtheorem{theo}{Theorem}

\newtheorem{subsecc}{}[secc]

                         \newtheorem{prop}{Proposition}

                         \newtheorem{Lemma}{Lemma}

                         \newtheorem{rema}{Remark}

                         \newtheorem{defin}{Definition}

                         \newtheorem{cor}{Corollary}

\def\a{\theta}
\newcommand{\T}{\mathbb{T}}

\newcommand{\R}{\mathbb{R}}
\newcommand{\Q}{\mathbb{Q}}
\newcommand{\Z}{\mathbb{Z}}
\newcommand{\N}{\mathbb{N}}

\def\carre{ \hfill $\Box$    }
\begin{document}
\title{Smooth linearization of commuting circle diffeomorphisms}
\author{Bassam Fayad and Kostantin Khanin}
\email{fayadb@math.univ-paris13.fr}
\email{khanin@math.toronto.edu}



\maketitle

\begin{abstract} We show that a finite number of commuting diffeomorphisms with simultaneously Diophantine  rotation numbers are smoothly conjugated to roations. This solves a problem raised by Moser in \cite{M}.\end{abstract}

\begin{secc} {\bf Introduction}  \rm

\medskip 

In this paper, we show that if a finite number of commuting smooth circle diffeomorphisms have  simultaneously Diophantine rotation numbers (arithmetic condition (\ref{dioph}) below), then the diffeomorphisms are smoothly (and simultaneously) conjugated to rotations (see Theorem \ref{theorem1} below).

The problem of smooth linearization of commuting circle diffeomorphisms was raised by Moser in \cite{M} in connection with the holonomy group of certain  foliations with codimension 1. Using the rapidly convergent Nash-Moser iteration scheme he proved that if the rotation numbers of the diffeomorphisms satisfy a simultaneous Diophantine condition and if the diffeomorphisms are in some $C^\infty$ neighborhood of the corresponding rotations (the neighborhood beign imposed by the constants appearing in the arithmetic condition, as usual in perturbative KAM theorems) then they are 
$C^\infty$-linearizable, that is, $C^\infty$-conjugated to rotations. 

In terms of small divisors, the latter result presented a new and striking phenomenon: if $d$ is the number of commuting diffeomorphisms, the rotation numbers of some or of all the diffeomorphisms may well be non-Diophantine,  but still, the full $\Z^d$-action is smoothly linearizable due to the absence of simultaneous resonances. Further, Moser shows in his paper that this new phenomenon is a {\it genuine} one in the sense that the problem cannot be reduced to that of a single diffeomorphism with a Diophantine frequency. indeed he shows that there exist numbers $\a_1,\ldots,\a_d$ that  are simultaneously Diophantine but such that for all linearly independent vectors  $a,b \in \Z^{d+1}$, the ratios  $(a_0+a_1\a_1+\ldots+ a_d\a_d) / (b_0+b_1\a_1+\ldots+b_d\a_d)$ are Liouville numbers. In this case, the theory for individual circle maps, even the global theorem of Herman and Yoccoz, does not suffice to conclude.

According to Moser, the problem of linearizing commuting circle diffeomorphisms could be regarded as a  model problem where KAM techniques can be applied to an overdtermined system (due to the commutation relations). This assertion could again be confirmed a quarter of a centuary later in a striking way by the recent work \cite{DK} where local rigidity of some higher rank abelian groups was established using a KAM scheme for an overdetermined system. 

At the time Moser was writing his paper, the {\it global} theory of circle diffeomorphisms was 
already known for a while, of which a highlight result is that a diffeomoprhism with a Diophantine rotation number is smoothly linearizable (without a {\it local} condition of closeness to a rotation, see \cite{Y}). The proof of the first global smooth linearization theorem given by Herman, as well as all the subsequent different proofs and generalizations, extensively used the Gauss algorithm of continued fractions that yields the best rational approximations for a real number.

As pointed out in Moser's paper, this is why the  analogue global problem for a commuting family of diffeomorphisms with rotation numbers satisfying a simultaneous Diophantine condition seemed difficult to tackle, due precisely to the absence of an analogue of the one dimensional continuous fractions algorithm in the case of simultaneous approximations of several numbers (by rationals with the same denominator).

Moser asked {\it under which conditions on the rotation numbers of $n$ smooth commuting circle diffeomorphisms can one assert the existence of a smooth invariant measure $\mu$? In particular is the  simultaneous Diophantine condition sufficient?}  Here, we answer this question positively (Theorem \ref{theorem1}, the existence of a smooth invariant measure being an equivalent statement to smooth conjugacy). On the other hand, it is not hard to see that the same arithmetic condition is optimal (even for the local problem) in the sense given by Remark \ref{liouv}.

Before we state our results and discuss the plan of the proofs,
we give a brief summary of the linearization theory of  single circle diffeomorphisms on which our proof relies.

We denote the circle by $\T =  \R / \Z$.  We denote by
${\rm Diff}^r_+(\T)$, $r \in [0, + \infty] \cup \lbrace \omega \rbrace$, the group of orientation
preserving
diffeomorphisms of the circle of class $C^r$ or real analytic. We represent the lifts  of these  diffeomorphisms as elements of $D^r(\T)$, the group of $C^r$-diffeomorphisms $\tilde{f}$ of the real line such that $f-{\rm Id}_\R$ is $\Z$-periodic. 

Following Poincar\'e, one can define the rotation number of a circle homeomorphism $f$ as the uniform limit $\rho_f = \lim (\tilde{f}^j(x)-x)/j {\rm mod } [1]$, where $\tilde{f}^j$ ($j \in \Z$) denote the iterates of a lift of $f$.   A rotation map of the circle with angle $\a$, that we denote by $R_\a: x \mapsto x + \a$, has clearly a rotation number equal to $\a$. Poincar\'e raised the problem of comparing the dynamics of a homeomorphism of the circle with rotation number $\a$ to the simple rotation $R_\a$. 

A classical result of Denjoy (1932) asserts that if $\rho_f=\theta$ is irrational (not in $\Q$) and if $f$ is of class $C^1$ and if $Df$ has bounded variations then $f$ is topologically conjugated to $R_\a$, i.e. there exists a circle homeomorphism $h$ such that $h \circ f \circ h^{-1} = R_{\theta}$. 

Considering the linearized version of the conjugation equation $H(x+\a) - H(x) = F(x)$ where $H$ and $F$ are real $\Z$-periodic functions  defined on $\R$ and where $F$ is assumed to have zero mean, it is easy to see (with Fourier analysis, due to the existence of the {\it small divisors}  
$|1-e^{i2\pi n \a}|$) that the existence of a smooth solution $H$, is guaranteed for all functions $F$ with zero mean if and only if $\a$ satisfies a Diophantine condition, i.e. if there exist $C>0$ and $\tau>0$ such that for any $k \in \Z$, $\|k\a\| \geq C {|k|}^{-\tau}$. Nonetheless, when $F$ is in some finite class of differentiability   and the linearized equation has a solution, this solution in general is of lower regularity than $F$. This is the so-called {\it loss of regularity} phenomenon. 

The first result asserting regularity of the conjugation of a circle diffeomorphism to a rotation was obtained by Arnol'd in the real analytic case:
if the rotation number of a real analytic diffeomorphism is Diophantine and if the diffeomorphism is sufficiently close to a rotation, then the conjugation is analytic.  The general idea, that is due to Kolmogorov, is to use a quadratic Newton approximation method to show that if we start with a map sufficiently close to the rotation it is possible to compose successive conjugations and get closer and closer to the rotation while the successive conjugating maps tend rapidly to the Identity. The Diophantine condition is used to insure that the loss of  differentiability in the linearized equation is fixed, which allows to compensate this loss at each step of the algorithm due to its quadratic convergence. Applying the same Newton scheme in the $C^\infty$ setting is essentially due to Moser.

At the same time, Arnol'd also gave examples of real analytic diffeomorphisms  with irrational rotation numbers for which the conjugating maps are not even absolutely continuous, thus showing that the {small divisors} effect was inherent to the regularity problem of the conjugation. Herman also showed that there exist "pathological" examples for any non-Diophantine irrational (i.e. Liouville) rotation number (see \cite[chap. XI]{H}, see also \cite{S}).

A crucial conjecture was that, to the contrary, the hypothesis of closeness to rotations should not be necessary for smooth linearization, that is, any smooth diffeomorphism of the circle with a Diophantine rotation number must be smoothly conjugated to a rotation. This {\it global} statement  was finally proved by Herman in \cite{H} for almost every rotation number, and later on by Yoccoz in \cite{Y} for all Diophantine numbers.

To solve the global conjecture, Herman, and later on Yoccoz, developped a powerful machinery giving sharp estimates on derivatives growth for the iterates of circle diffeomorphisms, the essential criterion for the $C^r$ regularity of the conjugation of a $C^k$ diffeomorphism $f$, $k \geq r \geq 1$, being the fact that the family of iterates $(f^n)$ should be bounded in the $C^r$ topology. The Herman-Yoccoz estimates on the growth of derivatives of the iterates of $f$ will be crucial for us in all the paper. 

\end{secc} 

\begin{secc} {\bf Results} \rm 
 
\medskip 

 For $\a \in  \T$ and  $r \in [1, + \infty] \cup \lbrace \omega \rbrace$,  we
denote by  ${\mathcal D}^r_\a$ the subset of  ${\rm Diff}^\infty_+(\T)$
of diffeomorphisms having rotation number $\a$. 

Let $d \in \N,$ $d \geq 2$, and assume that
$(\a_1,\ldots, \a_d) \in \T^d$ are such that there exist   $\nu>0$ and $C>0$ such
that for each  $k \in \Z^*$,
 \begin{eqnarray}\label{dioph}
{\rm max}(\|k\a_1\|,\ldots,\|k\a_d\|) \geq  C {|k |}^{-\nu}.
\end{eqnarray}

Finally, we say that a family of circle diffeomorphisms 
$(f_1,\ldots,f_d)$ is {\it commuting } if $f_i \circ f_j = f_j \circ
f_i$ for all $1\leq i \leq j \leq p$. Note that if $h$ is a homeomorphism of the circle such that $h \circ f_1
\circ h^{-1} = R_{\a_1}$, then for every $j \leq p$ we have that  $h \circ f_j
\circ h^{-1}$ commutes with $R_{\a_1}$, from which it is easy to see that  $h \circ f_j
\circ h^{-1} = R_{\a_j}$. Hence, for $r\geq 2$, Denjoy theory gives a homeomorphism that conjugates every $f_j$ to the corresponding rotation. Here, we prove the following.

\begin{theo}\label{theorem1} Assume that $\a_1,\ldots, \a_d$ satisfy
(\ref{dioph}) and let  $f_i  \in {\mathcal
D}_{\a_i}^\infty$, $i=1,\ldots,p$. If  $(f_1,\ldots,f_d)$ is commuting then, there exists  $h \in
{\rm Diff}^\infty_+(\T)$, such that for each $1\leq i \leq p$, $h \circ f_i
\circ h^{-1} = R_{\a_i}$. \end{theo}

\begin{rema} \rm \label{liouv}  Using Liouvillean constructions (constructions by successive conjugations) we see that the above sufficient arithmetic condition is also 
necessary to guarantee some regularity on the conujugating
homeomorphism $h$ (essentially unique, up to translation).  There is indeed a sharp dichotomy with the above statement in case the arithmetic condition is not satisfied (see for example \cite[chap. XI]{H}  and \cite{S} where the same techniques producing a single diffeomorphism readily apply to our context): {\it Assume that $\a_1,\ldots, \a_d$ do
  not satisfy
(\ref{dioph}), then there exist  $f_i  \in {\mathcal
D}_{\a_i}^\infty$, $i=1,\ldots,p$ such that  $(f_1,\ldots,f_d)$ is
commuting and such that the conjugating homeomorphism of the maps $f_i$ to the rotations $R_{\a_i}$ is not absolutely continuous.}
\end{rema}

As a corollary of Theorem \ref{theorem1} and of the local theorem (on commuting diffeomorphisms) of Moser in the real anlytic category \cite{M} we have by the same techniques as in \cite[chap. XI. 6]{H}:

\begin{cor}\label{analytic}  Assume that $\a_1,\ldots, \a_d$ satisfy
(\ref{dioph}) and let  $f_i  \in {\mathcal
D}_{\a_i}^\omega$, $i=1,\ldots,p$. If  $(f_1,\ldots,f_d)$ is commuting then, there exists  $h \in
{\rm Diff}^\omega_+(\T)$, such that for each $1\leq i \leq p$, $h \circ f_i
\circ h^{-1} = R_{\a_i}$. 
\end{cor}

If $(\a_1,\ldots, \a_d) \in \T^p$ are such that 
there exist $a \in (0,1)$ and infinitely many $k \in \N$ such that 
 \begin{eqnarray*} 
{\rm max}(\|k\a_1\|,\ldots,\|k\a_d\|) \leq a^k
\end{eqnarray*}
then it is possible to obtain with a construction by successive conjugations a commuting family
 $(f_1,\ldots,f_d) \in {\mathcal
D}_{\a_1}^\omega \times\ldots \times {\mathcal
D}_{\a_d}^\omega$ such that the conjugating homeomorphism of the maps $f_i$ to the rotations $R_{\a_i}$ is not absolutely continuous.

It is a delicate problem however to find the optimal arithmetic condition under which any commuting family of real analytic diffeomorphisms will be linearizable in the real analytic category. For a single real analytic diffeomorphism, the optimal condition was obtained by Yoccoz in \cite{Y3}.

\end{secc}

\medskip

\begin{secc} \rm {\bf Plan of the proof of Theorem \ref{theorem1}}

\medskip

 As in the global  theory of circle diffeomorphisms, we will start by proving the $C^1$ regularity of the conjugation and then we will derive from it by Hadamard convexity inequalities and bootstrap techniques the $C^\infty$ regularity. In each of these two moments of the proof the commutativity of the diffeomorphisms in question will be used differently. 

The first step in the proof is a simple arithmetic observation for which we need the following definition: given an angle $\a$ we say that a sequence of successive denominators of $\a$, $q_l,q_{l+1},\ldots,q_n$, is a {\it  Diophantine string of exponent} $\tau>0$ if for all $s\in [l,n-1]$, $q_{s+1} \leq q_{s}^\tau$. The observation is that if we consider a sufficiently large number of angles $\a_1,\ldots,\a_p$ such that each $d$-upple satisfies (\ref{dioph}) then we can  find Diophantine strings of the same exponent $\tau$ (function of $\nu$ and $d$) for different $\a_j$'s, such that these strings overlap (with a margin that can be made as large as the number of angles considered is large). In other words, one can follow successive denominators along a Diophantine string $i$ until its end, say at some $q_{j_i,n_i}$, where it is possible to switch to the next string $i+1$ starting from a denominator $q_{j_{i+1},l_{i+1}}$ that is well smaller than $q_{j_i,n_i}$ ($q_{j_{i+1},l_{i+1}} \leq q_{j_i,n_i}^{\xi}$, $\xi$ as small as desired as the number $p$ increases. The next elementary but crucial observation is that given $f_1,\ldots,f_d$ with rotation numbers   $\theta_1,\ldots,\theta_d$ satisfying (\ref{dioph}), it is possible, by considering compositions of these diffeomorphisms to obtain as much diffeomorphisms as desired with rotation numbers in such a way that any 
$d$-upple satisfies (\ref{dioph}). Sections \ref{sec4} and \ref{secalt} deal with these results on the alternated configuration of Diophantine strings. 

With this configuration in hand the proof of $C^1$-conjugacy goes as follows: first, to alleviate the notations we will consider only the case $d=2$ (the proof for $d \geq 3$ is exactly the same) and assume that the Diophantine strings of $\a=\rho_{f_1}$ and $\beta=\rho_{f_2}$ are themselves  in an alternated configuration (Conditions (\ref{11})--(\ref{bounds})) since this also does not make any difference with the proof in the general case.  If we denote by $m_n$ and $M_n$  the minimum and the maximum on the circle of $|x-f^{q_n}(x)|$ (where $q_n$ denotes the denominators of the convergents of $\a$, and with similar notations $\tilde{q}_n$,
$\tilde{m}_n$, and $\tilde{M}_n$ for $\beta$ and $g$), then a criterion for $C^1$-conjugacy of $f$ to a rotation is that $m_n/M_n$ be bounded. It is known that $m_n \leq \a_n \leq M_n$
where $\a_n =  |q_n \a - p_n|$ and the goal is to show that eventually both $m_n$ and $M_n$ become             comparable to $\a_n$ up to a multiplicative constant. In \cite{Y} a crucial recurrence relation between these quantities at the steps $n$ and $n+1$ is exhibited that allows to show, if a Diophantine condition holds on $\a$, that the quantities $m_n$ and $M_n$ end up having the same order. The latter recurrence relation is obtained as a result of the analysis of the growth of the Schwartzian derivatives of the iterates of $f$. 

Here we will rely on the same recurrence relation but use it only along the Diophantine strings and try to propagate the improvement of estimations when we switch strings using the commutation of $f$ and $g$. Actually this will work efficiently once we get started, namely once $M_s$ for $q_s$ in some Diophantine string for $\a$ is less than $\a_s^{1-\sigma}$ for some fixed $\sigma>0$ that depends on $\tau$ (it is possible to take $\sigma=1/(2\tau^2)$). This can be interpreted as a "local" result that yields $C^1$ conjugation for diffeomorphisms that are close to rotations (see Proposition \ref{local}). 
  
The existence of very long Diophantine strings (which corresponds to one of the angles being super-Liouvillean) presents the simplest case illustrating how the local situation can indeed be reached using only one string  (see Section \ref{special}). 

In general however, before reaching the local situation, switching from a string to a consecutive one may in fact lead to a worsening of the estimates or even to their complete loss (see the first equation in the proof of Lemma \ref{commutativity2}), so that a different strategy must be adopted. Keeping in mind that the objective is to see that $u_s \to 1$ where $u_s$ is such that $M_s=\a_s^{u_s}$ (with $\tilde{M}_s$, $\tilde{u_s}$, and $\beta_s$ for $\beta$ and $g$), the idea is to use each angle alone to study ``the dynamics'' of $u_s$: after we measure the gain in the exponent $u$ when we pass through a Diophantine string, we jump to the beginning of the successive string of the same angle. In this operation we can readily limit the loss in the exponent $u$ in function of the size of the jump (that in turn is less than the size of the overlapping Diophantine string of the other angle). Repeating these two steps inductively, we get a dynamics on the exponent $u_i$ measured at the exit of the $i^{\rm th}$ Diophantine string (of the same angle, see Lemma \ref{ll}). Doing so for each angle we see that at least for one of them, namely the one with the overall longest Diophantine strings (in the sense given by (\ref{A}) or (\ref{B})), the sequence $u_i$ (or $\tilde{u}_i$) eventually becomes larger than $1-\sigma$.  

The idea for proving higher regularity is to use convexity arguments as in \cite{H,Y} to bound the derivatives of the iterates of $f$ and $g$. To this difference that we will only seek to bound these derivatives for iterates  $f^u$ and $g^v$ at {\it Diophantine times} $u$ and $v$ that are (respectively) linear combinations of multiples of denominators $q_s$ and $\tilde{q}_s$ that belong to Diophantine strings (each $q_s$ is as usual multiplied by at most $q_{s+1}/q_s$). Due to the overlapping of strings, this will be  sufficient for proving regularity of the conjugation (see Section \ref{good}). 

Given a denominator $q_s$ in a Diophantine string,  the fact that the ratio
$q_{s+1}/q_s$ is bounded by a fixed power of
$q_s$ is naturally crucial in the control of the derviatives of the diffeomorphisms $f^{aq_s}, a\leq q_{s+1}/q_s$. Nonetheless, in the Herman-Yoccoz theory
for circle diffeomorphisms with Diophantine rotation number, the
control of the derivatives of $f^{q_s}$ itself are obtained using the
whole Diophantine condition on the diffeomorphism's rotation number (see 
the computations in \cite[section 8]{Y}). Still,
 we can see through the calculations of \cite[section 8]{Y} (see Section \ref{computations} below), that the existence of
a sufficiently long sequence of Diophantine string before and
up to some denominator $q_s$, combined with the existence of a $C^1$-conjugacy to a rotation, allows to give a bound on the derivatives
of $f^{q_s}$ that will be enough for our purpouse. 

Thus, in addition to the alternation of Diophantine strings used for $C^1$ we must make sure that there is enough Diophantine ``margin'' before $q_{l_i}$ and this is done (in Proposition \ref{conf}) through the use of even more numbers $\a_i$, which amounts to considering in the proof more diffeomorphisms of the form $f^i \circ g$. In a sense we use more and more relations in the commuting group of diffeomorphisms as we want to improve the regularity of the conjugation. 

The rest of the proof of higher regularity is inspired by the bootstrap calculations of \cite{Y}.

Nowhere in our proof of Theorem \ref{theorem1}, neither in the proof of the existence of $C^1$-conjugation nor in that  of its higher regularity, did we try to optimize on our {\it use of derivatives} of the diffeomorphism $f$, that is assumed to be of class $C^\infty$. For instance, the problem of finding the optimal regularity required on the diffeomorphisms that would guarantee $C^1$-conjugation under a given simultaneous Diophantine condition is an interesting problem that is not addressed in this paper.

\end{secc}

\medskip

\begin{secc}\label{sec4} \rm {\bf Preliminary : Diophantine strings}

\medskip

We recall that for every irrational number $\theta$ we can uniquely define an increasing  sequence of integers $q_n$ such that $q_1 =1$, and 
$$   \|  k \theta \|  >  \|  q_n  \theta\|, \quad \forall k<q_{n+1}, k \neq q_n.$$
This sequence is called the sequence of denominators of the best rational approximations, or convergents, of $\alpha$. 

Let $p \in \N$, and $\theta_1,\ldots, \theta_p$ be irrational
numbers. For $1 \leq j \leq p$, we denote by
$(q_{j,n})$ the sequence of denominators of the convergents of
$\theta_j$. For $\tau>0$, we define
$${\mathcal A}_{\tau}(\theta_j) = \lbrace s \in \N \  /  \ q_{j,s+1}
\leq q_{j,s}^\tau \rbrace.$$
A Diophantine string (with exponent $\tau$) for a number $\theta_i$
is then a sequence $l,l+1,\ldots,n-1  \in {\mathcal A}_{\tau}(\theta_i).$

We will prove in this section the main arithmetical result related with the simultaneous Diophantine
property (\ref{dioph}) that we will use to prove Theorem \ref{theorem1}.

\begin{prop}\label{mainarithmetic} Let $\nu>0$, $K>0$ and $d \in
  \N$, $d\geq 2$. There exists $p \in \N$ such that: if $\theta_1,\ldots,
  \theta_p$ are numbers for which there exists $C>0$ such that each
  $d-$upple (of disjoint numbers) 
  $(\theta_{i_1},\ldots,\theta_{i_d})$ satisfies (\ref{dioph}); if
  $ U>0$ is sufficiently large and if $U\leq V\leq U^K$; then there
  exists  $k \in \{1,\ldots,p\}$, with a Diophantine string
  $l,\ldots,n-1 \in {\mathcal A}_{\tau_{d-1}}(\theta_k)$ with 
$$q_{k,l} \leq U \leq V \leq q_{k,n},$$ 
where $(\tau_s)$ is the sequence defined by $\tau_0=\nu$ and
$\tau_s=2\tau_{s-1}+3$, for $s\geq 1$.
\end{prop}

\begin{defin}  For $\tau>0$, $C>0$, $d\in \N^*$, and an interval  $I \subset \R$ we define
$$D_{d,\tau,C}(I) = \{ (\theta_{1},\ldots,\theta_{d}) \in \R^d \ / \ {\sup_{1\leq i\leq d}} \| k \theta_i \| \geq C k^{-\tau},  \ \forall k \in I \}.$$
For $C=1$, we use the simplified notation $ D_{d,\tau}(I) := D_{d,\tau,1}(I). $
\end{defin}

We will need the following elementary but crucial arguing.

\begin{Lemma} Let $\nu>0, C>0, d\in \N, d\geq 2$. Define $\epsilon=1/(2\nu+2)$. There exists $U_0$ such that if $V \geq U \geq U_0$, and if $\theta_1,\ldots,\theta_d$ are numbers such that  
$$(\theta_1,\ldots,\theta_d) \in D_{d,\nu,C}([U,V]),$$
then if an integer $s \in [U,V]$ satisfies $\|  s \theta_d \|  \leq s^{-(2\nu+3)},$ we have that $$(\theta_{2},\ldots,\theta_{{d-1}}) \in D_{d-1,2\nu+3}([s,e])$$
with $e={\rm min} (V,{\|  s \theta_d \| }^{-\epsilon}).$ \label{el}
\end{Lemma}

\noindent{\it Proof.}  If $k \in [s,e]$ satisfies 
$${\sup_{i\leq d-1}} \|  k \theta_i\|  \leq k^{-(2\nu+3)},$$
then 
$${\sup_{i\leq d}} \|  ks \theta_i\|  \leq {(ks)}^{-(\nu+{1\over 2})},$$
which violates $(\theta_1,\ldots,\theta_d) \in D_{d,\nu,C}([U,V]),$
 if $s$ is sufficiently large.

\carre

Because $\eta=(2\nu+3)/(2\nu+2)>1$,  Lemma \ref{el} has the following immediate consequence.

\begin{cor} Let $\nu>0$, $K>0$ and $d \in
  \N$, $d\geq 2$. There exists $N \in \N$ such that: for each $C>0,$ there exists $U_0>0$, such that  if
  $ U \geq U_0$ and $U\leq V \leq U^K$, and if $p\geq N+d-1$ and $\theta_1,\ldots,
  \theta_p$ are numbers such that for each
  $d-$upple  (of disjoint indices) $i_1,\ldots,i_d$,
  $(\theta_{i_1},\ldots,\theta_{i_d}) \in D_{d,\nu,C}([U,V])$, then there exist $j_1,\ldots,j_{N} \leq p$ such that any $(d-1)$-upple (of disjoint numbers), $i_1,\ldots,i_{d-1} \in \{1,\ldots,p\}-\{j_1,\ldots,j_{N} \}$, satisfies $(\theta_{i_1},\ldots,\theta_{i_{d-1}}) \in D_{d-1,2\nu+3}([U,V])$.\label{corol}
\end{cor} 

\noindent{\it Proof.} We can in fact take $N= [ \ln K / \ln \eta]+2$. Let $p \geq N+d-1$ and let $k_1\in \N$, $k_1 \geq U$, be the smallest integer (if it exists) such that $\| k_1 \theta_i \|  \leq {k_1}^{-(2\nu+3)}$ for some $i \in \{1,\ldots,p\}$. Denote by $\theta_{j_1}$ the corresponding angle. 
Then, define $k_2 \geq  {\| k_1 \theta_i \|}^{-(2\nu+3)}$, to be the smallest integer (if it exists) such that  $\| k_2 \theta_i \|  \leq {k_1}^{-(2\nu+3)}$ for some $i \in \{1,\ldots,p\}-\{j_1\}$ and denote by $\theta_{j_2}$ the corresponding angle. Continuing this way, we construct a sequence ${j_1},\ldots,j_{N}$, and observe that $k_N \geq k_1^{\eta^N} >  V$. On the other hand,  Lemma \ref{el}  implies that any $(d-1)$-upple (of disjoint numbers), $i_1,\ldots,i_{d-1} \in \{1,\ldots,p\}-\{j_1,\ldots,j_{N} \}$, satisfies $(\theta_{i_1},\ldots,\theta_{i_{d-1}}) \in D_{d-1,2\nu+3}([k_s,  {\| k_s \theta_{j_s} \|}^{-(2\nu+3)}])$. But, by definition of $K_1,\ldots,k_N$, for every $i \in  \{1,\ldots,p\}-\{j_1,\ldots,j_{N} \}$, and for every $k \in  [U,k_1) \cup ({\| k_1 \theta_{j_1} \|}^{-(2\nu+3)}, k_{2}) \cup \ldots \cup  ({\| k_{N-1} \theta_{j_{N-1}} \|}^{-(2\nu+3)}, k_{N})$, we have $k \theta_i \geq {\| k \theta_{i} \|}^{-(2\nu+3)}$. Thus, $(\theta_{i_1},\ldots,\theta_{i_{d-1}}) \in D_{d-1,2\nu+3}([U,V])$. \carre

\bigskip

\noindent {\it Proof of Proposition \ref{mainarithmetic}.} If $p$ and $U$ are sufficiently large, applying Corollary \ref{corol} $d-1$ times (with $U^{1/(2\tau_{d-1})}$ instead of $U$) we get that there exists  $k \in 
\{1,\ldots,p\}$ such that $\theta_k \in D_{1,\tau_{d-1}}([U^{1/(2\tau_{d-1})},V])$. We claim that $\theta_k$ satisfies the properties required in Proposition \ref{mainarithmetic}. Indeed, it is sufficient to prove that $\theta_k$ must have a denominator $q_{k,l} \in [U^{1/(2\tau_{d-1})},U]$. But if this is not so, there is some $q_{k,l}\leq 
U^{1/(2\tau_{d-1})}$ such that $q_{k,l+1}\geq 
U$, but then $m=q_{k,l} U^{1/(2\tau_{d-1})} \leq U$ satisfies $\|  m \theta_k \| \leq m^{-\tau_{d-1}}$, in contradiction with  $\theta_k \in D_{1,\tau_{d-1}}([U^{1/(2\tau_{d-1})},V])$. \carre

\end{secc}

\medskip

\begin{secc}\label{secalt}  \rm {\bf Alternated configuration of denominators}

\medskip 

\begin{defin}\label{alter} We say that $\theta_1,\ldots, \theta_p$ are in an {\it
alternated configuartion} if there exist $\tau>1$, and two
increasing sequences of integers, $l_i$ and $n_i$ such that for each
$i$ there exists $j_i \in \lbrace 1,\ldots,p \rbrace$ with
  \begin{eqnarray}\label{2} l_{i},l_{i}+1,l_i+2,\ldots, n_{i}-1 \in
{\mathcal A}_{\tau}(\theta_{j_i}), \end{eqnarray}
and
\begin{eqnarray}\label{3} q_{{j_i,l_i}}^{\tau^2} \leq  q_{{j_i,n_i}}^{1 \over \tau^2} \leq {q}_{{j_{i+1},l_{i+1}}} \leq
q_{{j_i,n_i}}^{1 \over \tau}.  
\end{eqnarray}
\end{defin}

From Proposition \ref{mainarithmetic} it is straightforward to derive the following

\begin{prop}\label{conf} Let $\nu >0$, $\xi>0$, and $d \in
  \N$, $d\geq 2$. Let $\tau:= \tau_{d-1}$. There exists $p \in \N$ such that if  $\theta_1,\ldots,
  \theta_p$ are numbers for which there exists $C>0$ such that each
  $d-$upple (of disjoint numbers)
  $(\theta_{i_1},\ldots,\theta_{i_d})$ satisfies (\ref{dioph}) then
  $\theta_1,\ldots,
  \theta_p$ are in an alternated configuration (with
  exponent $\tau$) with   
in addition that for each $i$ there exists $l'_i$ such that 
 $q_{j_i,l'_i} \leq q_{j_i,l_i}^{\xi}$ and such that $l'_i,l'_i+1,\ldots,l_i-1 \in {\mathcal A }_{\tau}(\a_{j_i})$.  
\end{prop}

In our proof of Theorem \ref{theorem1}, we will show that if $f_1,\ldots,f_p$ are smooth commuting diffeomorphisms with rotation numbers $\theta_1,\ldots,\theta_p$ that are in an alternated configuration, then the diffeomorphisms are $C^1$-conjugated to rotations. The additional condition, i.e. the existence of long Diophantine strings before $q_{l_i}$ is then used to proof the higher regularity of the conjugacy, the higher the regularity required, the longer these Diophantine  strings should be ($\xi \rightarrow 0$). 

To adapt Proposition \ref{conf} to a family of $d$ commuting diffeomorphisms, we use the following somehow artificial trick\footnote{We may attribute, as we did in the introduction,  the usefulness of this trick to the fact that it exploits the raltions in the group, isomorphic to $\Z^d$, of commuting diffeomorphisms.}:  consider $\theta_1,\ldots,\theta_d$ satisfying (\ref{dioph})
and define for $s \in \N$ 
$$\tilde{\theta}_s = \theta_1+s\theta_2+\ldots+s^{d-1}\theta_d.$$
Observe that for any $p \geq d$, there exists $C>0$ such
that any disjoint indices $i_1,\ldots,i_d \leq p$, we have that $(\tilde{\theta}_{i_1},\ldots,
\tilde{\theta}_{i_d})$ satisfies (\ref{dioph}). Proposition \ref{conf} can
now be applied to $\tilde{\theta}_1,\ldots,\tilde{\theta}_p$. On the other hand, given $f_1,\ldots,f_d$ as in Theorem \ref{theorem1}, then the diffeomorphism 
$\tilde{f}_s={f_1} \circ f_2^s \circ f_3^{s^2} \circ \ldots \circ f_{d}^{s^{d-1}}$ has rotation number $\tilde{\theta}_s$.

Since it does not alter the proof but only alleviates the notations we
will assume for the sequel that $d=2$ and that $\a$ and $\beta$ are already in an
alternated configuration, that is,
there exist $\tau>1$,  and two increasing sequences of
integers, $l_i$ and $n_i$ such that
\begin{eqnarray}\label{11}
l_{2i},\ldots,n_{2i}-1 \in {\mathcal A}_{\tau}(\a) \\
l_{2i+1},\ldots,n_{2i+1}-1 \in {\mathcal A}_{\tau}(\beta)\label{22}
\end{eqnarray}
and

\begin{eqnarray}\label{bounds} q_{l_{2i}}^{\tau^2} \leq  q_{n_{2i}}^{1 \over \tau^2} \leq  \tilde{q}_{l_{2i+1}} \leq
q_{n_{2i}}^{1 \over \tau}, \quad  \tilde{q}_{l_{2i+1}}^{\tau^2} \leq  \tilde{q}_{n_{2i+1}}^{1 \over \tau^2} \leq  {q}_{l_{2i+2}} \leq
\tilde{q}_{n_{2i+1}}^{1 \over \tau}.  
\end{eqnarray}
where $(q_n)$ and $(\tilde{q}_n)$ denote respectively the sequences of denominators of the
convergents of $\a$ and $\beta$.

\end{secc}

\begin{secc} \rm {\bf Proof of $C^1$-conjugation}

\medskip

 Given $\a$
and $\beta$ satisfying (\ref{11})--(\ref{bounds})
 and two commuting diffeomorphisms $f \in  {\mathcal D}_\a$,  $g \in 
{\mathcal D}_\beta$ we will show in this section that $f$ and $g$ are $C^1$-conjugated
to the rotations $R_\a$ and $R_\beta$.

\begin{subsecc} \rm Let

\begin{eqnarray*}
\a_n &=& |q_n \a - p_n|, \ \ \ \ \ \ \ \ \ \ \ \ \ \ \ \ \beta_n =
|\tilde{q}_n \beta - \tilde{p}_n|  \\
M_{n} &=& \sup d(f^{q_{n}}(x),x),    \ \ \ \ \ \ \ \  \tilde{M}_{n} =
\sup d(g^{\tilde{q}_{n}}(x),x) \\
 m_{n} &=& \inf d(f^{q_{n}}(x),x), \ \ \  \ \ \  \   \  \ 
\tilde{m}_{n} = \inf d(g^{\tilde{q}_{n}}(x),x)
\\
U_n &=& {M_n \over m_n},  \ \ \ \ \ \ \ \  \ \ \ \ \ \ \  \ \ \ \ \ \
\  \  \ \tilde{U}_n = {\tilde{M}_n \over \tilde{m}_n} .
\end{eqnarray*}

Recall that  
\begin{eqnarray}
1/ (q_{n+1}+q_n) \leq \a_n \leq 1/q_{n+1}, \quad 
1/ (\tilde{q}_{n+1}+\tilde{q}_n) \leq \beta_n \leq 1/ \tilde{q}_{n+1}. \label{00}
\end{eqnarray}

Recall also that  since $\int_{\T}      |f^{q_n}-{\rm id} | d \mu = \a_n,$ (where $\mu$ is the unique probability measure invariant by $f$) then 
\begin{eqnarray*}
m_n \leq \a_n \leq M_n.
\end{eqnarray*}

Herman proved that a diffeomorphism is $C^r$ conjugated to a rotation
if and only if its iterates form a bounded sequence in the
$C^r$-topology (see \cite[Chap. IV] {H}). Based on the latter observation, the following
criterion for $C^1$ conjugacy was used in \cite{H} and in \cite[section 7.6]{Y}:

\begin{prop} If there exists $C>0$ such that $\limsup U_n \leq C$, then
$f$ is $C^1$-conjugated to $R_\a$ (actually $\liminf U_n \leq C$ is enough). \end{prop}

Our proof of $C^1$-conjugacy in Theorem  \ref{theorem1}  relies on the
following central estimate of \cite{Y}
\begin{prop}\label{Yoccoz} For any $f \in {\mathcal D}_{\a}$, for any
 $K\in \N$, there exists $C=C(f,K)$ such that
\begin{eqnarray}
M_n \leq M_{n-1} { (\a_n / \a_{n-1}) +CM_{n-1}^K \over 1-CM_{n-1}^{1/2}}\label{aaa} \\
m_n \geq m_{n-1} { (\a_n / \a_{n-1}) -CM_{n-1}^K \over 1+CM_{n-1}^{1/2}}.\label{bbb} 
\end{eqnarray}
\end{prop}
\end{subsecc}

\bigskip

\begin{subsecc} \rm {}  The goal of this section is to prove the following "local" result:

\begin{prop}\label{local} Let $\sigma = 1/(2\tau^2)$. There exists
  $i_0 \in \N$ such that if for some even (odd) integer  $i \geq i_0$,
  we have  
$$M_{n_i-1} \leq {1 \over q_{n_i}^{1-\sigma}}$$
(with $\tilde{M}_{n_i-1}$ and $ \tilde{q}_{n_i}$ instead of $M_{n_i-1}$ and  $q_{n_i}$ if $i$ is odd) then $U_n$ and
$\tilde{U}_n$ are bounded.
\end{prop}

\begin{rema} {\rm This can be viewed as a {local} result on
$C^1$-conjugation, since it states that if $M_{n_i-1}$ for $i$
sufficiently large is not too far from what it should be if $f$
were $C^1$-conjugated to the rotations, then $f$ and $g$  must indeed be
$C^1$-conjugated to the rotations.}
\end{rema}

\medskip

\noindent{\it Proof of Proposition \ref{local}.}   We will assume that
$i$ is even, the other case being similar. Due to the commutation of $f$ and $g$ we have

\begin{Lemma}\label{commutativity} Let $L_i = [\beta_{l_{i+1}-1} / \a_{n_i-1}]$, then
\begin{eqnarray} 
\label{01}
\tilde{M}_{l_{i+1}-1} &\leq& (1+L_i) M_{n_i-1} 
\\
\label{02}
\tilde{m}_{l_{i+1}-1} &\geq& L_i m_{n_i-1} 
\\
\label{03}
\tilde{U}_{l_{i+1}-1} &\leq& (1+{1\over L_i}) U_{n_i-1}. 
\end{eqnarray}
\end{Lemma}

\noindent {\sl Proof.} If we assume that 
${l_{i+1}}$ and ${n_i}$ (the other case being similar) we observe that for any $x \in \T$, 
\begin{eqnarray}\label{543} \bigcup_{k=0}^{L_i-1}R_\a^{kq_{n_i-1}}([x,R_\a^{q_{n_i-1}}(x)]) \subset 
[x,R_\beta^{\tilde{q}_{l_{i+1}-1}}(x)] \subset  
\bigcup_{k=0}^{L_i}R_\a^{kq_{n_i-1}}([x,R_\a^{q_{n_i-1}}(x)]).
\end{eqnarray}  
Since $f$ and $g$ commute there exists a continuous homeomorphism $h$ that conjugates $f$ to $R_\a$ and $g$ to $R_\beta$, and (\ref{01})--(\ref{03}) follow immediately from (\ref{543}). \carre

\medskip 

Proposition \ref{local} clearly follows from

\begin{Lemma}\label{commutativity2} Let $\sigma= 1/ {(2\tau^2)}$. There exists $i_0 \in \N$, such that if $i \geq i_0$ and $M_{n_i-1} \leq {1 / q_{n_i}^{1-\sigma}}$, then we have
\begin{eqnarray}\label{001} \tilde{M}_{n_{i+1}-1} \leq  {1 \over \tilde{q}_{n_{i+1}}^{1-\sigma}},
\end{eqnarray} 
and 
\begin{eqnarray}\label{002} \tilde{U}_{n_{i+1}-1} \leq  a_i U_{n_i-1} \end{eqnarray} 
with $a_i \geq 1$, and $\Pi_{i \geq i_0}  a_i < \infty$. 
\end{Lemma}

\noindent {\sl Proof of Lemma \ref{commutativity2}.} From (\ref{01}) we have 
\begin{eqnarray*} 
\tilde{M}_{l_{i+1}-1}  &\leq& (1+    {\beta_{l_{i+1}-1} \over \theta_{n_{i}-1}}) M_{n_i-1} \\
&\leq&  (1+2 {q_{n_i}  \over \tilde{q}_{l_{i+1}}})  {1 \over {{q}_{n_i}}^{1-\sigma}}   
\end{eqnarray*} 
hence  (\ref{bounds}) implies for $i$ sufficiently large 
\begin{eqnarray}\label{ml}  
\tilde{M}_{l_{i+1}-1}  \leq  {3 \over {\tilde{q}_{l_{i+1}}}^{1/2}}.   
\end{eqnarray} 

Now if we let $K= 2 [\tau]+2$ in Proposition \ref{Yoccoz}, then if $i \geq i_0$, $i_0$ sufficently large, we obtain from
(\ref{aaa}), (\ref{00}) and (\ref{22}) that 
\begin{eqnarray*} \tilde{M}_{l_{i+1}} &\leq& 
\tilde{M}_{l_{i+1}-1} {\beta_{l_{i+1}} \over \beta_{l_{i+1}-1}} (1+\tilde{q}_{l_{i+1}}^{-1/5})
\end{eqnarray*}
and by induction
\begin{eqnarray}\label{4r}  \tilde{M}_{n_{i+1}-1} \leq b_i  \tilde{M}_{l_{i+1}-1}  {\beta_{n_{i+1}-1} \over    \beta_{l_{i+1}-1}}  
\end{eqnarray}
with $b_i \geq 1$ and $\Pi_{i\geq i_0} b_i< \infty$. Thus, (\ref{001}) follows from (\ref{bounds}). 

By the same token, from (\ref{bbb}) in Proposition \ref{Yoccoz} and (\ref{ml}) we get  for $i \geq i_0$, $i_0$ sufficently large
\begin{eqnarray*}  \tilde{m}_{n_{i+1}-1} \geq c_i  \tilde{m}_{l_{i+1}-1}  {\beta_{n_{i+1}-1} \over    \beta_{l_{i+1}-1}}  
\end{eqnarray*}
with $c_i \leq 1$ and $\Pi_{i\geq i_0} c_i>0$. Together with (\ref{4r}) this implies that 
 \begin{eqnarray*}  \tilde{U}_{n_{i+1}-1} \leq d_i  \tilde{U}_{l_{i+1}-1} 
\end{eqnarray*}
with $d_i \geq 1$ and $\Pi d_i< \infty$. This, with (\ref{03}) and (\ref{bounds}), imply (\ref{002}). 
\carre

\bigskip

\end{subsecc}

\begin{subsecc}\label{special} \rm {\bf Moving towards the "local" situation. Proof of $C^1$-conjugation in a special case with long Diophantine strings.} The main ingredient in improving the bound of $M_{i}$ towardss the "local" condition of Proposition \ref{local} is the following. 

Let $A_i \geq \tau^4$ and $B_i  \geq \tau^4$ be such that 
$$q_{n_{2i}}= q_{l_{2i}}^{A_i}, \quad, \tilde{q}_{n_{2i+1}}= \tilde{q}_{l_{2i+1}}^{B_i}.$$

\begin{Lemma}\label{ll} For any $b \in \N$,  there exists $i_0$ such that  if  $i \geq i_0$ and $u_i>0$ is such that $M_{{l_{2i}} -1} = 1/ q_{l_{2i}}^{u_i}$, then we have 
$$M_{n_{2i}-1} \leq  1/ q_{n_{2i}}^{\rho_i}$$
with $\rho_i = {\rm min} (1-\sigma, A_i^b u_i)$.
\end{Lemma}

An immediate consequence of Proposition \ref{local} and Lemma \ref{ll} is the $C^1$-conjugacy in the particular case of very long Diophantine strings, namely if there exist $\epsilon>0$ and a strictly increasing subsequence 
 of the even integers ${(i_j)}_{\{j \in \N\}}$,   
such that
$$q_{n_{i_j}} \geq {q_{l_{i_j}}^{{(\|n q_{l_{i_j}})}^\epsilon}}.$$

\medskip

\noindent {\it Proof of Lemma \ref{ll}.} We denote $l=l_{i}$ and $n=n_{i}$. 
Let $r$ be such that 
$$ {A_i \over \tau^{4}} \leq  {\tau^{4r}} \leq  A_i.$$

Let $\tilde{K} := 2[\tau^{4b+1}]$, so that $\tilde{K}^r \geq A_i^b$. In Proposition \ref{Yoccoz} take $K:=[4\tau \tilde{K}]$. 

Notice that $q_{l}^{ {(\tau^{4r})}} \leq q_l^{A_i} \leq q_n.$ Hence, we can introduce a sequence of integers $p_s$, $s=0,\ldots,r,$ such that 
$p_0=l$, and for each $1\leq s \leq r$
$$q_{p_{s-1}}^{\tau^3} \leq q_{p_{s}} \leq q_{p_{s-1}}^{(\tau^4)}.$$

Using the first estimate of Proposition \ref{Yoccoz},
and following the idea of \cite[Sec. 7.4]{Y} it is easy to construct for
$j\in [l,n]$, positive sequences $u_j$ and $a_j \leq 2$ such that $u_l=1/\|n
q_l$, $a_l=1$, and for $j \in [l-1,n-1]$, $M_j \leq a_{j+1} /
q_{j+1}^{u_{j+1}}$, where for each $j \in [l,n-1]$ one of the two following alternatives holds:
\begin{itemize}
\item[{(i)}] If $\theta_{j}/\theta_{j-1} \leq CM_{j-1}^{K/2}$ then $a_{j+1}=a_j$ and $u_{j+1}=\tilde{K} u_j$;
   
\item[{(ii)}] If $\theta_{j}/\theta_{j-1} > CM_{j-1}^{K/2}$ then $a_{j+1}=b_ja_j$ and $u_{j+1}= u_j$, with $\Pi b_j \leq 2$. In this case, we actually have $M_{j} \leq b_j M_{j-1} \a_{j} / \a_{j-1}$.   
\end{itemize}

Now, if there exists $s \in [0,r-1]$ such that for every $j\in [p_s,p_{s+1}-1]$, alternative (ii) holds, then 
 (assuming without loss of generality that $\tau \geq 2$) we have
\begin{eqnarray*} M_{p_{s+1}-1} &\leq& 2 M_{p_s-1} { \a_{p_{s+1}-1} \over   \a_{p_s-1} } \\ 
 &\leq&   {q_{p_s} \over  q_{p_{s+1}}} \\
 &\leq&   {1 \over  q_{p_{s+1}}^{1-\sigma}}
\end{eqnarray*}
after which, and as in the proof of Lemma \ref{commutativity2},  only
alternative (ii) can happen for all $j \in [p_{s+1}-1,n-1]$, so that,  arguing again as in Lemma \ref{commutativity2}, we get $M_{n-1} \leq {1 \over q_{n}^{1-\sigma}}$ and we finish. 

Otherwise, we have for every $s \in [0,r-1]$, at least one  $j\in [p_s,p_{s+1}-1]$ for which alternative (i) holds, hence $u_{p_{s+1}} \geq \tilde{K} u_{p_s}$. Subsequently,  $u_{p_{r}} \geq \tilde{K}^r u_{l} \geq A_i^b u_l$.   The Lemma is thus proved. \carre 

\end{subsecc}

\medskip 

\begin{subsecc} \rm {\bf Proof  of $C^1$-conjugation in the general case.} Recall that $A_i \geq \tau^4$ and $B_i  \geq \tau^4$ are such that 
$$q_{n_{2i}}= q_{l_{2i}}^{A_i}, \quad, \tilde{q}_{n_{2i+1}}= \tilde{q}_{l_{2i+1}}^{B_i}.$$

Then, clearly at least one of the following two limits holds 
\begin{eqnarray}  \label{A}
\limsup {\Pi_{j=1}^i { A_j^2} \over \Pi_{j=1}^i { B_j}} = +\infty \\
 \label{B}
\limsup {\Pi_{j=1}^i { B_j^2} \over \Pi_{j=1}^i { A_j}} = +\infty 
\end{eqnarray}

We will assume that (\ref{A}) holds, the other case being similar. 
We will show how  Lemma \ref{ll} applied with $b=2$, implies that eventually the condition of Proposition \ref{local} will be satisfied, thus yielding $C^1$-conjugacy.

Notice first that $q_{l_{2(i+1)}} \leq q_{n_{2i}}^{B_i}$. Furthermore, $M_{{l_{2(i+1)}} -1} \leq M_{n_{2i}-1}$ since $q_{l_{2(i+1)}} \geq {q}_{n_{2i}}^{\tau}$. 

Now, if $i_0$ is some sufficienlty large ineteger, and if  at step $i_0$ we do not have $M_{n_{2i_0}-1} \leq  1/ q_{n_{2i_0}}^{1-\sigma}$, we observe as above that  $M_{{l_{2(i_0+1)}} -1}  \leq M_{n_{2i_0}-1} \leq 1 / q_{l_{2(i_0+1)}}^{u_{i_0} A_{i_0}^2/B_{i_0}} $. A continued application of the Lemma hence shows that either at some  $i\geq i_0+1$  the condition of 
Proposition \ref{local} will be satisfied, or for every $i\geq i_0$, $M_{{l_{2i}} -1}  \leq 1 / q_{l_{2i}}^{u_{i_0} \Pi_{j=i_0}^{i-1} (A_{j}^2/B_{j})} $ which, with our assumption that (\ref{A}) holds,  contradict the fact that for every $i$, 
$M_{{l_{2i}} -1}  \geq 1 / (2q_{l_{2i}})$.

\medskip 

\begin{rema} \rm In the general situation, the alternated configuration of denominators may require the use of more than two angles, that is more than two diffeomorphisms. Our proof remains quite the same. Indeed, let $\a_1,\ldots,\a_p$ be in an alternated  configuration as in definition \ref{alter}. Define $A_i$ such that  $q_{{j_i,n_i}} = q_{{j_i,l_i}}^{A_i}$. Then there exists $k \in [1,p]$ such that 
$$\limsup_{I \in \N} {\Pi_{j_i=k, i\leq I} A_i^{p+1} \over \Pi_{j_i\neq k, i\leq I} A_i} = + \infty$$
and the proof of $C^1$-conjugation follows the same lines as above with this difference that we would take $b=p+1$ in Lemma \ref{ll},  $\tilde{K}= 2[\tau^{4b+1}]$) and then $K:=[4\tau \tilde{K}]$ in Proposition \ref{Yoccoz} which is possible since the diffeomorphisms we are considering are of class $C^\infty$. We see here the dramatic increase in our need of differentiability to prove $C^1$-conjugation as the number $d$ of           commuting diffeomorphisms in Theorem \ref{theorem1} increases.
\end{rema}

\end{subsecc}

\end{secc}

\begin{secc} \rm {\bf Higher  regularity}

\medskip

We fix $r \geq 2$. Knowing that the diffeomorphisms $f$ and $g$ are $C^1$-conjugated to the rotations, we will now prove that the conjugacy is in fact  of class $C^r$. 

In all the sequel, we fix  $k=[(r+2)(2+\tau)]+2$. And we take $\xi = 1/k$ in Proposition \ref{conf}.

As in the proof of $C^1$-conjugation, we will continue to assume for simplicity that  we are given $\a$
and $\beta$ satisfying (\ref{11})--(\ref{bounds}) with in addition that there exists for each $i$, $l'_i$ such that 
if $i$ is even, then
\begin{eqnarray}\label{subs} q_{l'_i} \leq q_{l_i}^{1/k}, \quad {\rm and } \  l'_i,\ldots,l_i-1 \in {\mathcal A }_\tau(\a),
\end{eqnarray}
 with a similar property involving $\beta$ if $i$ is odd.

Given two commuting diffeomorphisms $f \in  {\mathcal D}_\a$,  $g \in 
{\mathcal D}_\beta$ such that $f$ and $g$ are $C^1$-conjugated, we will show that the conjugacy is actually of class $C^r$.

\begin{subsecc} \label{good} \rm {\bf The control of the derivatives at alternating  "Diophantine times" is sufficient.}  We define two sets of integers, the "Diophantine times", as

\begin{eqnarray*} {{\mathcal A}} &=& \lbrace m \in \N \ / \ m= \sum a q_s  \ {\rm with } \ s \in [l_{2i},n_{2i}-1], \ i \in \N,  \ a \leq  q_{s+1} / q_s \rbrace \\ 
  \tilde{{\mathcal A}} &=& \lbrace m \in \N \ / \ m= \sum a \tilde{q}_s  \ {\rm with } \ s \in [l_{2i+1},n_{2i+1}-1], \ i \in \N,  \ a \leq  \tilde{q}_{s+1} / \tilde{q}_s \rbrace. \end{eqnarray*}

We also define two sets of diffeomorphisms 
\begin{eqnarray*} 
{\mathcal Z} &=& \lbrace f^n \ / \  n \in \N \rbrace \\
{\mathcal C} &=& \lbrace f^u \circ g^v \ / \  u \in {\mathcal A}, v \in \tilde{\mathcal A} \rbrace. 
\end{eqnarray*}

The following is an elementary Lemma due to (\ref{bounds})

\begin{Lemma}\label{un} If we denote
$${\mathcal O} = \lbrace u \a + v \beta {\rm mod}[1] \ / \ u \in {{\mathcal A}}, v \in \tilde{{\mathcal A}} \rbrace$$
then $\overline{\mathcal O} = \T$.

As a consequence, we have that ${\mathcal C}$ is dense in ${\mathcal Z}$ in the $C^0$-topology.

\end{Lemma}

It follows from the above Lemma that it is enough to control the derivatives of the $f^u$ and $g^v$ at the {\it  Diophantine times} $u \in {\mathcal A}$ and $v \in \tilde{{\mathcal A}}$: 

\begin{cor}\label{cr} If  ${\mathcal C}$ is bounded in the $C^{r+1}$-topology, then the conjugating diffeomorphism $h$ of $f$ to $R_\a$ is of class $C^r$. 
\end{cor}

\noindent {\it Proof.} We know that ${1 \over n} \sum_{i=0}^{n-1} f^i$
converges in the $C^0$-topology to $h$ (see \cite[chap. IV]{H}). From Lemma \ref{un}, this implies that there exist sequences $(u_n)$ and $(v_n)$ of numbers in ${{\mathcal A}}$ and  $\tilde{{\mathcal A}}$ such that  the sequence ${1 \over n} \sum_{i=0}^{n-1} f^{u_i} \circ g^{v_i} $ converges in the $C^0$ topology to $h$. By our $C^{r+1}$-boundness assumption, we can extract from the latter sequence a sequence that converges in the $C^r$-topology, so that necessarily $h \in {\rm Diff}^r_+(\T)$. \carre

\end{subsecc}

\bigskip

\begin{subsecc}\label{computations} \rm  It follows from standard computations (see  \cite[section 8.10]{Y}) that the assumption of  corollary \ref{cr} holds true if we prove 

\begin{Lemma}\label{deux}  There exists $\nu>0$ such that, for $i$ (even) sufficenlty large, we have for any $s \in [l_i,n_i-1]$ and for any $0 \leq a \leq q_{s+1} / q_s$
$${\| \ln D f^{aq_s} \|}_{r+1} \leq q_s^{-\nu}$$
(with $g$ and $ \tilde{q}_{s}$ instead of $f$ and  $q_{s}$ if $i$ is odd).
\end{Lemma}

\noindent {\it Proof. }
We will only work with $f$ since the arguments for $g$ are the same. The proof is based on the estimates of \cite[section 8]{Y} and  we start by recalling some facts that were proven there:

For $k \in \N^*$,  define for $s \in \N$, $\Delta_{s}^{(k)} = \| D^{k-1} \ln Df^{q_s} \|_0 + \a_s$. Then it follows from the $C^1$-conjugation of $f$ to $R_\theta$  (see \cite[lemme 5]{Y}) that 
$$\Delta_{s}^{(k)} \leq q_s^{(k-1)/2}.$$
We will use this fact with   $k=[(r+2)(2+\tau)]+2$ and use the notation $\Delta_{s} $  for $\Delta_{s}^{(k)} $.

Observe that  for $s \in [l'_i,n_i-1]$, we have if $i$ is sufficiently large 
\begin{eqnarray}\label{six}  {(\Delta_{s} q_{s+1})}^{1/k} q_s^{-1} \leq q_s^{-1/4}. \end{eqnarray}

Hence it follows from \cite[lemme 14 in section 8.8]{Y} that  for any $s\in [l'_i,n_i-1]$, 
 and for any $0 \leq a \leq q_{s+1} / q_s$, we have 
\begin{eqnarray}\label{aq} {\| \ln D f^{aq_s} \|}_{r+1} \leq C q_s^{-1} {(\Delta_{s} q_{s+1})}^\rho \end{eqnarray} 
where $ \rho = (r+2) / k$ and $C$ is some constant. 

If we denote
$$\Delta'_{s} = {\rm Sup} \lbrace {| (D^{k-1} \|n D f^{q_t} \circ f^m) {(Df^m)}^{k-1} |}_0, \ 0 \leq t \leq s, m \geq 0 \rbrace,$$
then we have $\Delta_s \leq C \Delta'_s$ for some constant
$C$. Observe that since ${\|Df^m\|}_0$ is bounded we have $\Delta'_{s}
\leq C q_s^{(k-1)/2}$ for some constant $C$. If we denote $V_s = {\rm
  Max} \lbrace \Delta'_t / q_t, \ 0 \leq t \leq s \rbrace$, then due
to (\ref{six}) we have from \cite[section 8.9]{Y} that for $s \in   [l'_i,n_i-1]$, 
$$V_{s+1} \leq V_s (1+ C q_s^{-1/4})$$
for some constant $C$. Hence, for $s \in [l'_i,n_i-1]$, we have $V_s \leq 2 V_{l'_i} \leq C q_{l'_i}^{(k-3)/2}$. If $s \geq l_i$ this gives 
\begin{eqnarray*}
\Delta'_{s} &\leq& C q_{s} q_{l'_i}^{(k-3)/2}  \\
&\leq& C q_{s}^2  
\end{eqnarray*}
because we assumed that $q_{l'_i} \leq q_{l_i}^{2/k}$.

Finally, if $s \in [l_i,n_i-1]$ we have that  ${(\Delta_{s}
  q_{s+1})}^\rho \leq C q_s^{(2+\tau)(r+2)/k} \leq q_s^{1-1/k}$ and we
conclude using (\ref{aq}) that the statement of Lemma \ref{deux} holds
which ends the proof of higher regularity. \carre

\end{subsecc}

\end{secc}

\medskip 

\noindent {\sc Acknowledgments.} The first author is grateful to Arthur Avila, Rapha\"el Krikorian, and Jean-Christophe Yoccoz for important suggestions and comments. The second author is grateful to Alexey Teplinsky for many useful conversations.

\medskip

\frenchspacing
\bibliographystyle{plain}

\medskip

\end{document}